# Reduction of Look Up Tables for Computation of Reciprocal of Square Roots

Shadrokh Samavi and M.Reza Jahangir

Department of Electrical and Computer Engineering, Isfahan University of Technology, Isfahan 84156-83111, Iran.

*Abstract*— Among many existing algorithms, convergence methods are the most popular means of computing square root and the reciprocal of square root of numbers. An initial approximation is required in these methods. Look up tables (LUT) are employed to produce the initial approximation. In this paper a number of methods are suggested to reduce the size of the look up tables. The precision of the initial approximation plays an important role in the quality of the final result. There are constraints for the use of a LUT in terms of its size and its access time. Therefore, the optimization of the LUTs must be done in a way to minimize hardware while offering acceptable convergence speed and exactitude.

## I. INTRODUCTION

The reciprocal of square root is an important operation for applications, such as graphics or scientific computations. Most multimedia applications and three-dimensional graphics require reciprocal of square root. For many computational processes, including vector normalization [6], least squares lattice filters [10], Cholesky decomposition [2], and Givens rotations [5], a square root is first computed and then used as the divisor in a subsequent divide operation. A more efficient method for performing this computation is to first compute the inverse square root and then use it as the multiplier in a subsequent multiply operation [6]. Users of most processors till a decade ago had to use software methods to compute reciprocal of square roots [12], [17]. Software methods do not yield high quality and real time frame production. Because of its usefulness in 3D graphics applications, special instructions for inverse square root have been added to a number of microprocessors, such as MIPS R10000/R12000 [9], [11] and IBM PowerPC740/750 [4]. This instruction is implemented by a square rooting followed by a division using a divide/square-root unit based on digit-recurrence methods, or by a convergence method such as Newton-Raphson (NR) method using a multiplier [6],[7],[13]. A convergence algorithm such as Newton-Raphson by utilizing look up tables gives better cost-performance results compared to a digit-recurrence method [16].

Originally, convergence methods were known as functional iteration algorithms. They became popular due to their small computational delay [14]. A typical convergence process is shown in Figure 1.

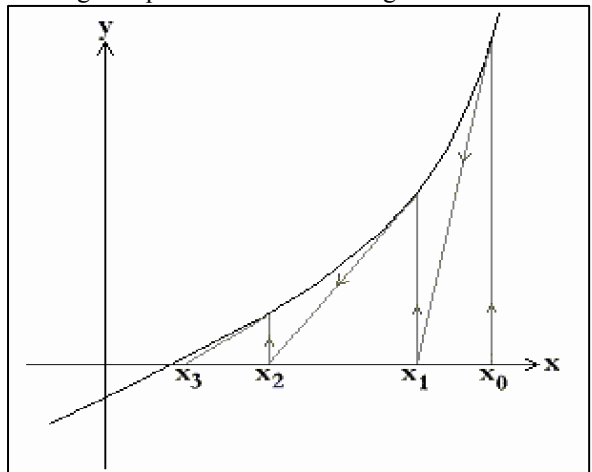

Figure 1. Convergence process in Newton-Raphson algorithm.

The iterative equation for the Newton-Raphson method is:

$$X_{n+1} = X_n - \left[\frac{f(X_n)}{f'(X_n)}\right] \qquad (1)$$

where $X_0$ is the initial approximation and $X_n$ approaches the root of the function $f(x)$.

Convergence methods start with an initial approximation and use a number of iterations to converge toward the final result. Usually, iterative methods have quadratic convergence [3] [14]. This means that the number of valid digits of the result doubles after every iteration. On the other hand, digit-by-digit methods produce one digit of the result at every iteration. Therefore, they are said to have linear convergence.

Some of the basic characteristics of the convergence algorithms include 1) no need for division operation, 2) quadratic convergence order, 3) maximum error of less than one "unit in least position" (ULP) and 4) use of look up tables. To obtain the required precision the important factors in the implementation of a convergence algorithm are the number of bits of the initial approximation and the number of iterations. Therefore, either a higher number of bits must be stored in look up tables or higher number of iterations must be performed. For example, in an IEEE single precision number, where a 23-bit mantissa is used, 6-bit look up table content and 2 iterations are necessary.

For a double precision case, 7 bits in the content area of the look up table and three iterations are required. In average, at every iteration step, two add/sub operations and four multiplications are required.

During the past decades, the use of ROM based LUTs in the computational circuits, microprocessors and math coprocessors have increased dramatically. This is due to high device density of the ROMs, more regular layout structure of the memories, and lower implementation and test costs of memories compared to other devices [6], [12].

One of the cardinal factors for a better performance of a convergence method is the appropriate use of LUTs for generation of initial approximation. An appropriate initial point will result in reduced number of iterations and more precise result. Figure 2 shows the quadratic convergence of the Newton-Raphson algorithm implemented using a 64×6 LUT.

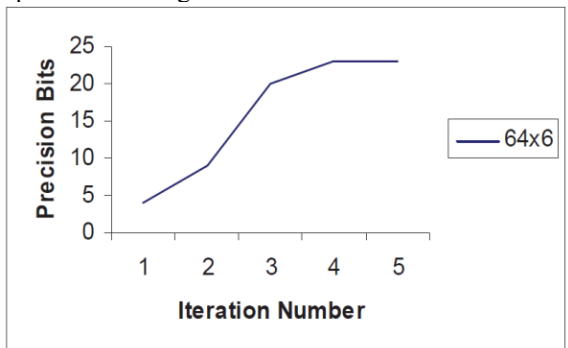

Figure 2. Quadratic convergence of NR algorithm.

It can be deduced from Figure 2 that after every iteration of the NR algorithm, the number of valid result bits is doubled. After three iterations the initial 4 bits of precision reaches the desired 23 bits.

In the following sections, a number of suggestions are offered to improve the implementation of the algorithms and their simulation results are presented.

## II. OPTIMIZED LOOK UP TABLE

Since there are constraints in use of large memory spaces in ICs, it is imperative to optimize the size of these tables while keeping the required precision and performance. There are two approaches to achieve this goal. The first is to omit the fixed bits that are always 1 or 0. Also bits that have a regular pattern can be omitted. The second approach is to find points using interpolation within spans that the function has linear characteristics. Therefore, a smaller number of points can be directly stored in the table. The first approach reduces the length of each entry while the second approach reduces the length of the whole table.

One of the methods to generate the initial approximation is to find the product of a table entry with some modified form of the bit pattern of the operand [17]. These are called auxiliary tables because they do not directly contain the reciprocal of square root. Keeping in perspective parameters such as error in the first iteration, amount of divergence and average number of iterations to get the required precision level, the simulation results suggest that an auxiliary table of 4096×25 size is the smallest and optimized one. Table 1 summarized the simulation results of this part and shows the reason for this selection.

Table 1. Simulation results for various auxiliary tables with 25-bit contents

| # of address bits | 6 | 8 | 10 | 12 | 14 | 16 |
|---|---|---|---|---|---|---|
| Error | Acceptable | Acceptable | Acceptable | Acceptable | Acceptable | Acceptable |
| Divergence | 0.23% | 0.31% | 0.02% | None | None | None |
| Average # of iterations | 1.006 | 1.009 | 1.001 | 1.0 | 1.0 | 1.0 |

The meaning of "acceptable" in Table 1 refers to an error that is less than 1 ULP. The reason for choosing 25-bit words for the content of the LUT is its better convergence and lower errors in the first or second iteration steps of the NR algorithm. It is obvious that with each bit increase in the number of address bits, the size of the LUT is doubled. To reduce this size we suggest the use of linear interpolation method. This has never been used for the auxiliary tables.

## III. INTERPOLATION IN AUXILIARY TABLES

The goal of interpolation is to estimate the value of a function between two known points. The estimation error must be small enough to be disregarded. Replacing the function in a small range with a straight line is the interpolation. The process has low precision and is only suggested to be used in spans that the function has linear behavior [9], [11].

In storing the values of f(x) in a LUT, instead of placing all of the points in the table, only every other point is stored. In order to compute the function at a point $t$, when $t$ is not directly in the table, two points on two sides of $t$ are retrieved from the table. Then $f(t)$ is computed using the following equation:

$$t_i \leq t \leq t_{i+1}$$
$$f(t) = f(t_i) + (t - t_i) \cdot \left[\frac{f(t_{i+1}) - f(t_i)}{t_{i+1} - t_i}\right] \quad (2)$$

The goal of the performed simulations has been to come up with an optimized size of LUT and an appropriate interpolation factor. This factor refers to the distance between two consecutive stored values. An appropriate interpolation will result in an error that is less than one ULP while the number of diverged results and that of iteration steps are minimized.

In the followings we look at the effects of different parameters on factors such as error, divergence and number of iterations.

An auxiliary table with 25-bit content and address sizes of 12, 14 and 16 bits were simulated for 10,000 random numbers with IEEE 754 format. Also, different interpolation factors were used and for each case the number of diverged entries is found. Figure 3 shows simulation results for the NR algorithm.

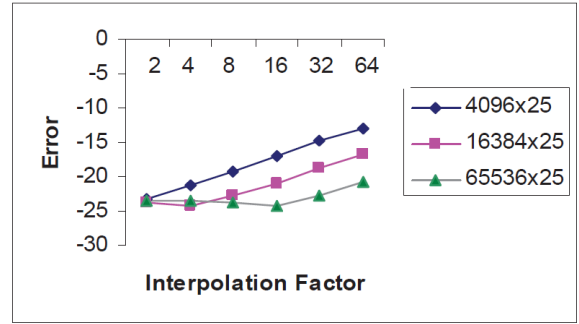

Figure 4. Average error as function of interpolation factor in NR algorithm.

It can be noted that as the number of address bits of the table is increased the average error of the first iteration is reduced. For the 64K table the interpolation factors of 2, 4, 8 and 16 generate acceptable error. The interpolation factor of 2 generates results that are independent of the size of the tables. This is the common point of all three graphs of Figure 4.

The required number of iterations to produce an acceptable error in the NR algorithm is presented in Figure 5.

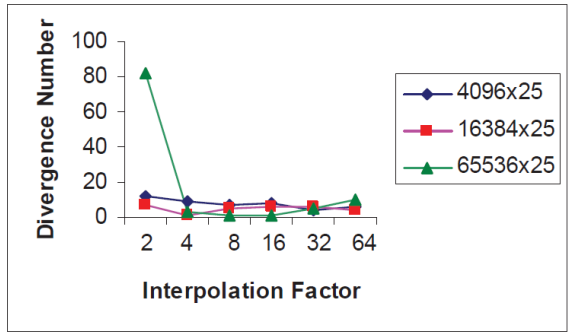

Figure 3. The effect of interpolation factor on divergence.

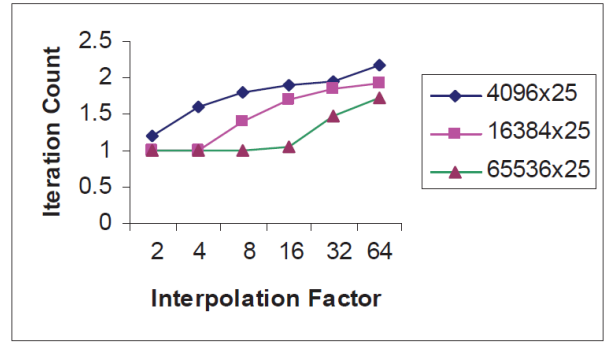

Figure 5. Required number of iterations for different tables.

The simulation results show that except for interpolation factors of 2 and 64 the divergence of a 64K×25 is much better than the other tables and is close to zero. For the 4K×25 we see a better divergence at 2 and 64 interpolation factors. Since the worst-case divergence is 0.82% it can be deduced that the divergence is not dependent on the interpolation factor.

The tolerable error is one that is less than one ULP. This, for IEEE 754 single precision standard, is an error of less than $2^{-23}$. In the following plots we show the error with its exponent part, for example, an error of $2^{-23}$ is represented by –23.

For 10,000 samples the effect of linear interpolation factor on the error, after one iteration, for different auxiliary tables were simulated. The results of the simulation are presented in Figure 4.

The simulation results specify that the larger the table the faster is the convergence of the algorithm. There is a relationship between the amount of error caused by interpolation factor and the average number of iterations. From Figures 4 and 5 it is apparent that for the 64k table the average error and average iteration follow the same pattern.

The simulation results for the auxiliary table are summarized in Table 2. It can be concluded that an auxiliary table of 4k×25 with an interpolation factor of 2 is optimum in terms of error, number of iterations and divergence.

It can also be concluded that an interpolation factor of 2 can reduce the size of this table to half while the error will be kept below one ULP and the number of iterations and diverged points are not dramatically increased.

Table 2. Characteristics of a 4k×25 auxiliary table.

| Interpolation factor | 2 | 4 | 8 | 16 | 32 | 64 |
|---|---|---|---|---|---|---|
| Error after 1st iteration | Acceptable | Unacceptable | Unacceptable | Unacceptable | Unacceptable | Unacceptable |
| Error after 2nd iteration | Acceptable | Acceptable | Acceptable | Acceptable | Acceptable | Acceptable |
| Divergence | 0.12% | 0.09% | 0.07% | 0.08% | 0.04% | 0.06% |
| # of iterations | 1.20 | 1.60 | 1.80 | 1.90 | 1.95 | 2.16 |

## IV. USE OF INTERPOLATION IN MAIN TABLES

Most of implementations of reciprocal of square root have used auxiliary tables to come up with the initial approximation. One of the major drawbacks of this method is its need for a complex multiplication. This multiplication is done to retrieve the word corresponding to the input bits. An approach for improving convergence methods is the use of what we call the "main lookup tables" (MLT) as opposed to "auxiliary" lookup tables (ALT). In that case the need for the costly multiplication is diminished. The simulation results show that MLTs are more efficient than ALTs.

Among the advantages of MLTs one can mention no occurrence of divergence in the final results, more accuracy and also smaller number of iterations. It can be deduced from the simulation results that a MLT with 11 bits of address and 23-bit words would perform optimally. The sheer size of this table shows an advantage over a best size 4095×25 ALT. The corresponding ALT had more than twice the capacity of the suggested MLT.

The content of MLT linearly grows therefore interpolation can be used to reduce the size of these tables. This has no precedence in the literature. An appropriate interpolation is one that does not exceed the error beyond one ULP, increase the number of divergent points or increase the number of iterations.

In the following section the simulation results are presented. The effect of parameters such as number of address bits on parameters such as error level, divergence and number of iterations are studied.

A number of simulation runs were performed to see the effect of the size of MLT on the number of divergent points in the Newton-Raphson method. MLTs with 11, 12 and 16 bit address lines and 23-bit words were simulated for different interpolation factors. Figure 6 shows the simulation results.

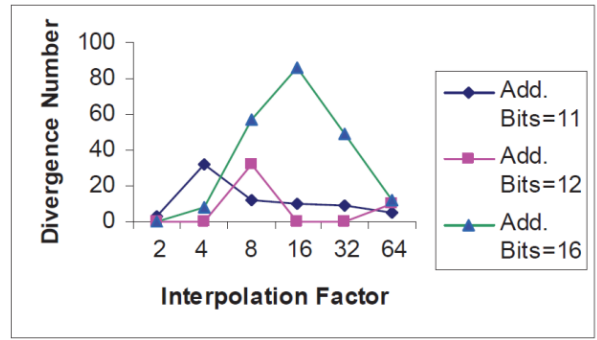

Figure 6. Effect of the interpolation factor on divergence of the algorithm.

It can be seen that the three MLTs follow the same patterns. Divergence peaks at a certain address length and reduces afterwards.

Maximum divergence occurs at a MLT with 16-bit address with interpolation factor of 16. The maximum divergence for 11-bit table is at interpolation factor of 4 and for 12-bit table at interpolation factor of 8. Divergence of zero occurred for a 12-bit table at interpolation factors of 2, 4, 16 and 32.

The effect of interpolation factor on the error level is the next subject of investigation. MLTs with 11, 12 and 16-bit addressing and 23-bit words were used. Figure 7 shows the result of using different interpolation factors. A notable point is that all simulated cases resulted in errors less than one ULP. The minimum error occurred for the table with 12-bit address and maximum error belonged to MLT with 16-bit addressing. The effect of the interpolation factor on the average error is very much similar to the effect of the interpolation on the divergence. In another words the minimum error is at the interpolation factor of 2 and then this error increases linearly to a peak where it goes on a decline afterwards.

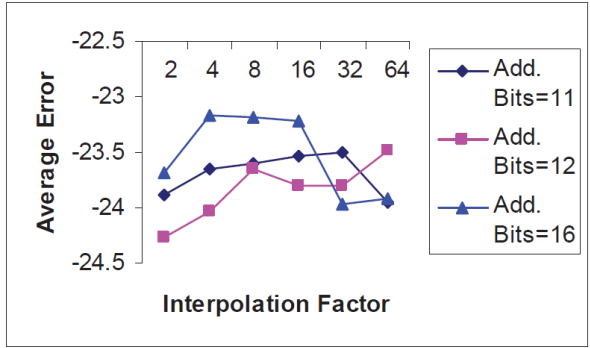

Figure 7. The effect of the interpolation factor on the average error.

Another parameter of interest is the number of iterations and its variations as a function of the interpolation factor. Simulations were performed on MLTs with 11, 12 and 16 address lines and 23-bit words. Figure 8 shows the results of the simulation.

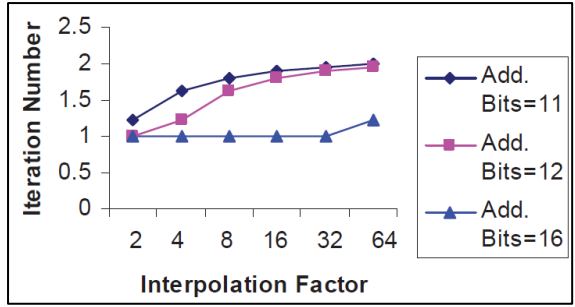

Figure 8. The effect of interpolation factor on the number of iterations.

It appears that as the number of address bits increases the required number of iterations decreases. The effect of interpolation factor is different for different tables. In the table with 16 address lines the number of iterations is fixed at 1 up to the interpolation factor of 32 and beyond that the speed of convergence is reduced. For the other tables the reduction in the convergence speed is linear as the interpolation factor increases.

It can be concluded from the simulation results of this section that two of the MLTs can perform optimally. The table with 12 addresses lines and interpolation factor of 4 (MLT A) as well as the table with 11 address lines and interpolation factor of 2 (MLT B) performed better than the others. Both tables are of the same physical size. MLT A creates the least number of divergent points while the advantage of MLT B is because of the least number of iterations.

Table 3 shows the interpolation characteristics of these two MLTs. The implementation of a MLT with interpolation factor of 2 is easier than that of 4, therefore, it is suggested to use MLT B. Hence it is concluded that by using interpolation the size of an optimum MLT can be reduced by 50 percent without increasing the error beyond one ULP or noticeable increase in the number of iterations.

Table 3. Characteristics of MLTs A and B

| Interpolation factor | 2 | 4 | 8 | 16 | 32 | 64 |
|---|---|---|---|---|---|---|
| Final error for MLT A | Acceptable | Acceptable | Acceptable | Acceptable | Acceptable | Acceptable |
| Number of iterations for MLT A | 1.0 | 1.235 | 1.614 | 1.810 | 1.904 | 1.954 |
| Final error for MLT B | Acceptable | Acceptable | Acceptable | Acceptable | Acceptable | Acceptable |
| Number of iterations for MLT B | 1.216 | 1.614 | 1.807 | 1.902 | 1.954 | 1.997 |
| Divergence for MLT B | 0.03% | 0.32% | 0.12% | 0.1% | 0.09% | 0.05% |

## V. WORD LENGTH REDUCTION

Besides reducing the length of a LUT by means of interpolation the width of the table can be reduced. Eliminating unnecessary bits does this. These bits can be categorized in two groups. Bits that are constant form the first category. In the second category there are bits that have a regular pattern. Joint application of both methods would reduce both length and width of a LUT. In the following, two examples for reducing the width of a table are presented.

The IEEE standard of 754 for floating point numbers requires 23-bit mantissa. Therefore, the word size of the MLT is 23 bits long. The initial approximation of the square root reciprocal has 3 most significant bits that can be eliminated.

The first MSB bit is always "1" and can be removed. The second MSB follows a regular pattern. For a 2K word MLT up to address 1593 and for a 4K word MLT up to address 3186 all of the second MSB bits are "1" and thereafter they are zero. Therefore, by comparing the most significant bits of the input with a fixed combination the mentioned bit can be constructed. This eliminates the need for storing this bit.

The third MSB of the content of every line of a MLT follows a regular pattern too. This bit in a 2K word MLT up to address 627 and in a 4K word MLTs up to address 1254 is equal to "1" and thereafter it is the complement of the second MSB. Hence, by using most significant bits of the address both second and third bits of the content of any MLT can be constructed.

VI. CONCLUSION

In this paper parameters that are affected by the size of LUTs such as computational error and number of iterations were analyzed. Newton-Raphson algorithm as an example of convergence algorithm was used in these analysis but the results could be generalized to all convergence algorithms. The simulation results suggested that a 2K-word table would perform optimally. It was also shown that interpolation could reduce the size of a table without jeopardizing the precision of the results. Finally, further reduction of the tables was suggested by trimming the words.